\documentclass[a4paper]{article}

\pdfoutput=1 

\usepackage[T1]{fontenc}
\usepackage[utf8]{inputenc}
\usepackage[english]{babel}
\usepackage{graphics}
\usepackage{graphicx}
\usepackage{epsfig}
\usepackage{amsmath}
\usepackage{amssymb}
\usepackage{amsthm}
\usepackage{url}
\usepackage{hyperref}
\usepackage{float}
\usepackage{pifont}
\usepackage[table]{xcolor}
\usepackage{lmodern}
\usepackage{xspace}
\usepackage{bbold}
\usepackage[capitalize,nameinlink]{cleveref}
\usepackage{subcaption}
\usepackage{mathtools}
\usepackage{adjustbox}

\title{GR decompositions and their relations to Cholesky-like factorizations}

\author{
Peter Benner\thanks{%
Computational Methods in Systems and Control Theory,
Max Planck Institute for \newline Dynamics of Complex Technical Systems,
Sandtorstr.~1, 39106 Magdeburg, Germany
}
\and Carolin Penke$^*$ \thanks{Corresponding author: \url{penke@mpi-magdeburg.mpg.de}}}

\usepackage{times,cite}
\usepackage[T1]{fontenc}
\usepackage[utf8]{inputenc}
\usepackage{times,cite}
\usepackage[T1]{fontenc}
\usepackage[utf8]{inputenc}
\usepackage{amssymb}
\usepackage{algorithm}
\usepackage{algorithmicx}
\usepackage{algpseudocode}
\usepackage{pgfplots}
\usepackage{tikz}
\usepackage{caption}
\usepackage{subcaption}
\usepackage{mymacros}
\usepackage{matrix}
\usepackage{enumitem}
\usepackage{changepage}

\setlist[enumerate]{topsep=0pt,itemsep=-1ex,partopsep=1ex,parsep=1ex}

\usepackage{graphicx}

\begin{document}

\maketitle

\paragraph{Abstract} ~\\

%
%
For a given matrix, we are interested in computing GR decompositions $A=GR$, where $G$ is an isometry with respect to given scalar products. The orthogonal QR decomposition is the representative for the Euclidian scalar product. For a signature matrix, a respective factorization is given as the hyperbolic QR decomposition. Considering a skew-symmetric matrix leads to the symplectic QR decomposition. The standard approach for computing GR decompositions is based on the successive elimination of subdiagonal matrix entries. For the hyperbolic and symplectic case, this approach does in general not lead to a satisfying numerical accuracy. An alternative approach computes the QR decomposition via a Cholesky factorization, but also has bad stability. It is improved by repeating the procedure a second time. In the same way, the hyperbolic and the symplectic QR decomposition are related to the $LDL^T$ and a skew-symmetric Cholesky-like factorization. We show that methods exploiting this connection can provide better numerical stability than elimination-based approaches.

\section{Introduction and Preliminaries}\label{intro}
Bilinear forms on $\mathbb{R}^{m}$ with respect to a nonsingular matrix $M\in\mathbb{R}^{m\times m}$ are defined as $\langle x, y \rangle_M= x^T My$ \cite{HigMMetal05}.
The adjoint of a matrix $A\in\mathbb{R}^{m\times m}$ is given as $A^{\star_M}\in\mathbb{R}^{m\times m}$ and is uniquely defined by $\langle Ax, y \rangle_M = \langle x, A^{\star_M}y \rangle_M$ for all $x,y\in\mathbb{R}^{m}$. The adjoint generalizes the transpose $.^{\tran}$. It holds $A^{\star_M} = M^{-1}A^{\tran}M$. Similarly, the adjoint of a rectangular matrix $A\in\mathbb{R}^{m\times n}$ is given with respect to two bilinear forms induced by matrices $M\in\mathbb{R}^{m\times m}$ and $N\in\mathbb{R}^{n\times n}$ as $A^{\star_{M,N}}$ \cite{HigMT10}. It is defined by satisfying the identity $\langle Ax, y \rangle_M = \langle x, A^{\star_{M,N}}y \rangle_N$ and it holds $A^{\star_{M,N}} = N^{-1}A^{\tran}M$. We are interested in computing decompositions $A=GR\in\mathbb{R}^{m\times n}$, $m\geq n$, where $G\in\mathbb{R}^{m\times n}$ is an $(M,N)$-isometry, i.e. $G^{\star_{M,N}}G=I_n$ and $R\in\mathbb{R}^{n\times n}$ \cite{Wat07}. $(M,N)$-isometries are useful for devising structure-preserving methods, for example in the context of eigenvalue computations \cite{Bun86}. This work considers bilinear forms in real space but the theory is easily extended to complex space or sesquilinear forms. The most well known representative of this class of decompositions is the (thin) QR decomposition. Here $M=I_m$, $N=I_n$ and $R$ is upper triangular. With respect to these matrices, an isometry is a matrix with orthonormal columns. Typically, the QR decomposition is computed in a stable fashion by successively eliminating subdiagonal entries of the matrix using orthogonal transformations. The decomposition has a well known connection to the Cholesky factorization. Let $A$ have full column rank. It holds that $A=QR$ is a thin QR decomposition if and only if $R$ defines a Cholesky decomposition $R^TR = A^TA$. Computing $Q:=AR^{-1}$ provides an alternative to the column elimination approach. For tall and skinny matrices, this method has a much lower computational effort but is known to be unstable. However, the stability can be drastically improved by doing a second repetition, i.e. compute the QR decomposition of $Q$ \cite{YamNYF15}. This is also done in the context of $\mathcal{H}$-matrices
\cite{Lin04}. In this work, we investigate whether this observation also holds for other QR-like decompositions.

\section{GR decompositions and Cholesky-like factorizations}
We now consider a scalar product induced by a signature matrix $
  \Sigma_m = \diag{\sigma_1,\dots,\sigma_m}$, where $\sigma_1,\dots,\sigma_m \in \{+1,-1\}$.  A $(\Sigma_m,\Sigma_n)$-isometry $H$ is called hyperbolic and fulfills the property 
  $H^{\tran}\Sigma_m H = \Sigma_n$. For a given $\Sigma_m$, the hyperbolic QR decomposition  $A=HR$, where $H\in\mathbb{R}^{m\times n}$, $R\in\mathbb{R}^{n\times n}$ upper triangular, exists if all principal submatrices of $A$ are nonsingular \cite{Wat07}. It can be computed via successive column elimination, similar to the orthogonal case \cite{Hou15}. The diagonal values of $\Sigma_n$ are determined by the used transformations and are a subset of the diagonal values of $\Sigma_m$. The role of the Cholesky factorization is now played by the $LDL^{\tran}$ factorization. $A=HR$ is a hyperbolic QR decomposition with respect to $\Sigma_m$ and $\Sigma_n$ if and only if $R^{\tran}\Sigma_n R = A^{\tran}\Sigma_m A$ gives a scaled $LDL^{\tran}$ factorization. As the computation of the $LDL^{\tran}$ factorization can be unstable, one typically relies on the slightly altered Bunch–Kaufman  factorization \cite{AshGL99}. Here, $D$ is allowed to have $2\times 2$ diagonal blocks and pivoting is introduced in form of a permutation $P$. Using this factorization as a starting point, we arrive at a different variant of the $HR$ decomposition, called indefinite QR decomposition in \cite{Sin06}.  It can be computed from the Bunch-Kaufman (BK) factorization in the following way.
 
 \begin{enumerate}
  \item Compute BK factorization $A^{\tran}\Sigma A =: PLDL^{\tran}P^{\tran}$.
  \qquad 2. Diagonalize $D=:V\Lambda V^T$,\ $\Sigma_n := \sign{\Lambda}$
  \item[3.] $R := |\Lambda|^{\frac{1}{2}}V^{\tran}L^{\tran}P^{\tran}$, $H := AR^{-1}$
 \end{enumerate}
 In the resulting QR-like decomposition, $R$ is no longer upper triangular, but block upper-triangular with permuted columns. In Figure \ref{Fig:NumRes} we see how applying this method a second time to $H$ affects the numerical accuracy.

 We now consider scalar products induced by $J_m= \begin{bmatrix}
            0 & I_m\\
            -I_m & 0
           \end{bmatrix}$. A $(J_m,J_n)$-isometry $S\in\mathbb{R}^{2m\times 2n}$ fulfills the property $S^{\tran}J_mS = J_n$ and  is called symplectic. The symplectic QR decomposition $A=SR$ can again be computed by the successive introduction of zeros in the columns using symplectic transformations. In contrast to the hyperbolic and orthogonal QR decomposition, $R$ is not upper triangular but of the form $R=\left[
   \begin{matrix}
             \resizebox{\widthof{$R_{11}$}}{!}{$\utri$} & \resizebox{\widthof{$R_{11}$}}{!}{$\smbsutri{0}$} \\
             \resizebox{\widthof{$R_{11}$}}{!}{$\smbsutri{0}$}  & \resizebox{\widthof{$R_{11}$}}{!}{$\utri$}\\
            \end{matrix}\right]=P_s^{\tran}\hat{R}P_s$, where $\hat{R}$ is upper triangular with $2\times 2$ diagonal blocks on the diagonal and  $P_s=[e_1,e_3,\dots,e_{2n-1},e_2,e_4\dots,e_{2n}]$ is the perfect shuffle. 
            This variant of the symplectic QR decomposition corresponds to the skew-symmetric Cholesky factorization described in \cite{BenBFetal00,Bun82}. $A=SR$ is a symplectic QR decomposition if and only if $\hat{R}^{\tran}(P_sJ_nP_s^{\tran}) \hat{R}$ is a Cholesky-like decomposition of the skew-symmetric matrix $P_sA^{\tran}J_m AP_s^{\tran}$. Similar to the hyperbolic case, pivoting in form of a permutation matrix $P$ can be introduced to increase the stability of the Cholesky-like decomposition, leading to an altered SR decomposition $AP^{\tran}P_s=SR$. A Cholesky-based computation takes the form
            
\begin{enumerate}
  \item Compute Cholesky-like factorization $A^{\tran} J_m A =: P\hat{R}^{\tran} P_sJ_nP_s^{\tran} \hat{R}P $.
  \quad 2. $R:=P_s^{\tran} \hat{R}P_s$. \quad3. $S:= AP^{\tran}P_SR^{-1}$.
 \end{enumerate}
Again there exists the possibility to repeat the procedure for the computed symplectic factor $S$.
            
   \newlength\figureheight
\newlength\figurewidth 
\setlength\figureheight{0.2\textwidth}
\setlength\figurewidth{\textwidth}
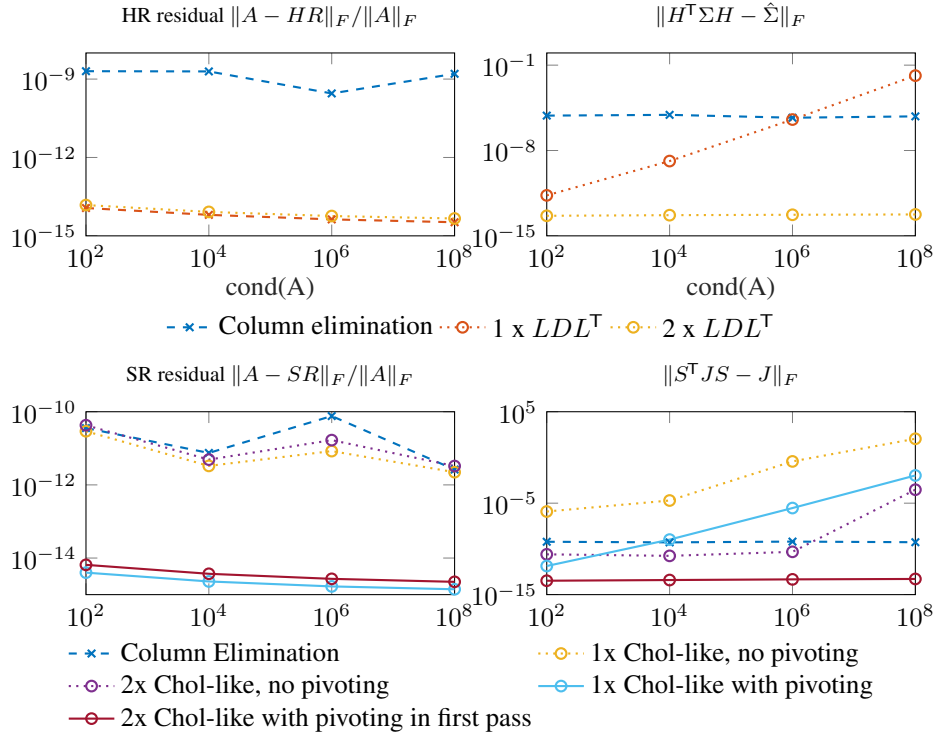
\begin{figure}[H]
%
%
\definecolor{mycolor1}{rgb}{0.00000,0.44700,0.74100}%
\definecolor{mycolor2}{rgb}{0.85000,0.32500,0.09800}%
\definecolor{mycolor3}{rgb}{0.92900,0.69400,0.12500}%
\definecolor{mycolor4}{rgb}{0.49400,0.18400,0.55600}%
\definecolor{mycolor5}{rgb}{0.46600,0.67400,0.18800}%
\definecolor{mycolor6}{rgb}{0.30100,0.74500,0.93300}%
\definecolor{mycolor7}{rgb}{0.63500,0.07800,0.18400}%

\begin{tikzpicture}
\pgfplotsset{compat=1.3}
\pgfplotsset{every axis/.append style={xlabel shift = -3pt }}

\begin{axis}[%
ylabel near ticks,
every axis plot/.append style={thick},
width=0.4\figurewidth,
height=\figureheight,
at={(0\figurewidth,0\figureheight)},
scale only axis,
xmode=log,
xmin=100,
xmax=100000000,
xtick={100, 10000,1000000,100000000},
xlabel style={font=\color{white!15!black}},
xlabel={cond(A)},
ymode=log,
ymin=1e-15,
ymax=1e-8,
yminorticks=true,
ylabel style={font=\color{white!15!black}},
axis background/.style={fill=white},
title={\footnotesize HR residual $\|A-HR\|_F/\|A\|_{F}$}
]

\addplot [color=mycolor1, dashed, mark=x, mark options={solid, mycolor1}, forget plot]
  table[row sep=crcr]{%
100	2.00355771464019e-09\\
10000	1.96205933255133e-09\\
1000000	2.82335592217485e-10\\
100000000	1.59239079151354e-09\\
};
\addplot [color=mycolor2, dashed, mark=x, mark options={solid, mycolor2}, forget plot]
  table[row sep=crcr]{%
100	1.163715322965e-14\\
10000	6.3598378214208e-15\\
1000000	4.26472421829671e-15\\
100000000	3.31607946750606e-15\\
};
\addplot [color=mycolor3, dotted, mark=o, mark options={solid, mycolor3}, forget plot]
  table[row sep=crcr]{%
100	1.50702463840337e-14\\
10000	8.24600979866679e-15\\
1000000	5.73777261704727e-15\\
100000000	4.62118255113158e-15\\
};

\end{axis}

\begin{axis}[%
ylabel near ticks,
width=0.4\figurewidth,
height=\figureheight,
at={(0.5\figurewidth,0\figureheight)},
scale only axis,
xmode=log,
xmin=100,
xmax=100000000,
xminorticks=true,
xtick={100, 10000,1000000,100000000},
xlabel style={font=\color{white!15!black}},
xlabel={cond(A)},
ymode=log,
ymin=1e-15,
ymax=1,
yminorticks=true,
axis background/.style={fill=white},
title={\footnotesize$\|H^{\tran}\Sigma H-\hat{\Sigma}\|_F$},
legend style={at={(-.2,-.5)}, anchor=center, legend cell align=left, align=left, draw=white!15!black,  legend columns=3, draw=none},
every axis plot/.append style={thick},
]
\addplot [color=mycolor1, dashed, mark=x, mark options={solid, mycolor1}]
  table[row sep=crcr]{%
100	7.26302377466114e-06\\
10000	8.38340432806101e-06\\
1000000	4.87575264375706e-06\\
100000000	6.29763200724347e-06\\
};
\addlegendentry{Column elimination}

\addplot [color=mycolor2, dotted, mark=o, mark options={solid, mycolor2}]
  table[row sep=crcr]{%
100	2.08787670877151e-12\\
10000	1.35980738609837e-09\\
1000000	3.44976465021385e-06\\
100000000	0.0137736157546972\\
};
\addlegendentry{1 x $LDL^{\tran}$}

\addplot [color=mycolor3, dotted, mark=o, mark options={solid, mycolor3}]
  table[row sep=crcr]{%
100	4.41545065859879e-14\\
10000	4.9634125844834e-14\\
1000000	5.34662729414506e-14\\
100000000	5.69452647053831e-14\\
};
\addlegendentry{2 x $LDL^{\tran}$}

\end{axis}

\end{tikzpicture}%
%
%
\definecolor{mycolor1}{rgb}{0.00000,0.44700,0.74100}%
\definecolor{mycolor2}{rgb}{0.85000,0.32500,0.09800}%
\definecolor{mycolor3}{rgb}{0.92900,0.69400,0.12500}%
\definecolor{mycolor4}{rgb}{0.49400,0.18400,0.55600}%
\definecolor{mycolor5}{rgb}{0.46600,0.67400,0.18800}%
\definecolor{mycolor6}{rgb}{0.30100,0.74500,0.93300}%
\definecolor{mycolor7}{rgb}{0.63500,0.07800,0.18400}%

\begin{tikzpicture}
\pgfplotsset{compat=1.3}
\pgfplotsset{every axis/.append style={tick label style={font=},xlabel shift = -3pt }}

\begin{axis}[%
ylabel near ticks,
every axis plot/.append style={thick},
width=0.4\figurewidth,
height=\figureheight,
at={(0\figurewidth,0\figureheight)},
scale only axis,
xmode=log,
xmin=100,
xmax=100000000,
xtick={100, 10000,1000000,100000000},
xlabel style={font=\color{white!15!black}},
xlabel={cond(A)},
ymode=log,
ymin=1e-15,
ymax=1e-10,
yminorticks=true,
ylabel style={font=\color{white!15!black}},
axis background/.style={fill=white},
title={\footnotesize SR residual $\|A-SR\|_F/\|A\|_{F}$}
]
\addplot [color=mycolor1, dashed, mark=x, mark options={solid, mycolor1}, forget plot]
  table[row sep=crcr]{%
100	3.5703988265145e-11\\
10000	7.46861962534135e-12\\
1000000	7.62649333813563e-11\\
100000000	2.62682736790492e-12\\
};
\addplot [color=mycolor3, dotted, mark=o, mark options={solid, mycolor3}, forget plot]
  table[row sep=crcr]{%
100	2.92789711217255e-11\\
10000	3.30331151396089e-12\\
1000000	8.39531463428923e-12\\
100000000	2.22376690302738e-12\\
};
\addplot [color=mycolor4, dotted, mark=o, mark options={solid, mycolor4}, forget plot]
  table[row sep=crcr]{%
100	4.27399969236871e-11\\
10000	4.91046164377375e-12\\
1000000	1.67450158937541e-11\\
100000000	3.27611439097022e-12\\
};
\addplot [color=mycolor6, mark=o, mark options={solid, mycolor6}, forget plot]
  table[row sep=crcr]{%
100	3.99230716123302e-15\\
10000	2.2919314194965e-15\\
1000000	1.67736588708107e-15\\
100000000	1.40078475057222e-15\\
};
\addplot [color=mycolor7, mark=o, mark options={solid, mycolor7}, forget plot]
  table[row sep=crcr]{%
100	6.51593778585303e-15\\
10000	3.71213348504926e-15\\
1000000	2.70050143763618e-15\\
100000000	2.23697897691211e-15\\
};
\end{axis}

\begin{axis}[%
ylabel near ticks,
width=0.4\figurewidth,
height=\figureheight,
at={(0.5\figurewidth,0\figureheight)},
scale only axis,
xmode=log,
xmin=100,
xmax=100000000,
xminorticks=true,
xtick={100, 10000,1000000,100000000},
xlabel style={font=\color{white!15!black}},
xlabel={cond(A)},
ymode=log,
ymin=1e-15,
ymax=100000,
yminorticks=true,
ylabel style={font=\color{white!15!black}},
axis background/.style={fill=white},
title={\footnotesize $\|S^{\tran}JS-J\|_F$},
legend style={at={(-.2,-.5)}, anchor=center, legend cell align=left, align=left, draw=white!15!black, font=, legend columns=2, draw=none},
every axis plot/.append style={thick},
ylabel shift = -3 pt
]
\addplot [color=mycolor1, dashed, mark=x, mark options={solid, mycolor1}]
  table[row sep=crcr]{%
100	6.02544665872194e-10\\
10000	5.17351389690564e-10\\
1000000	6.29363795433229e-10\\
100000000	5.38011563329762e-10\\
};
\addlegendentry{Column Elimination}


\addplot [color=mycolor3, dotted, mark=o, mark options={solid, mycolor3}]
  table[row sep=crcr]{%
100	1.22232399042358e-06\\
10000	1.95782898588639e-05\\
1000000	0.369047634965527\\
100000000	109.971052085946\\
};
\addlegendentry{1x Chol-like, no pivoting}

\addplot [color=mycolor4, dotted, mark=o, mark options={solid, mycolor4}]
  table[row sep=crcr]{%
100	2.54346848785739e-11\\
10000	1.73244740285142e-11\\
1000000	4.79169637550855e-11\\
100000000	0.000297392237709787\\
};
\addlegendentry{2x Chol-like, no pivoting}


\addplot [color=mycolor6, mark=o, mark options={solid, mycolor6}]
  table[row sep=crcr]{%
100	1.34303358276298e-12\\
10000	1.04443527121516e-09\\
1000000	2.93659549270909e-06\\
100000000	0.0107321854888783\\
};
\addlegendentry{1x Chol-like with pivoting}

\addplot [color=mycolor7, mark=o, mark options={solid, mycolor7}]
  table[row sep=crcr]{%
100	3.3069350627552e-14\\
10000	3.93877650167981e-14\\
1000000	4.65185292244139e-14\\
100000000	5.16324647057781e-14\\
};
\addlegendentry{2x Chol-like with pivoting in first pass}


\end{axis}

\end{tikzpicture}%
 \caption{Numerical results for computing decompositions of randomly generated matrices of size $500\times 500$ (HR) or $1000\times 1000$ (SR).} \label{Fig:NumRes}
 \end{figure} 
 Figure \ref{Fig:NumRes} shows how the accuracy of HR and SR decompositions can be improved by computing them via the $LDL^{\tran}$ and the Cholesky-like decomposition. Computing the decompositions via factorizations employing pivoting with two iterations leads to accuracies that do not derail for badly conditioned matrices.

\vspace{\baselineskip}


\end{document}